\documentclass[12pt,reqno]{amsart}

\theoremstyle{definition}

\theoremstyle{plain}

\theoremstyle{remark}

\begin{document}

\begin{center}
{\sc T.I. Gaisin}\\
\smallskip
{\large \bf The maximum principle for manifolds
over a local algebra (2).}
\end{center}
\bigskip
\par
Let $A$ be a finite-dimensional local commutative algebra over
$R$, $\dim_RA=n$. In this work we consider compact manifolds over
$A$, and prove that the real part of an $A$-differentiable
function is constant. Also we find estimates for the dimensions of
some spaces of 1-form.

For the theory of manifolds over algebras we refer the reader to
[1],[2],[3]. We assume that all manifolds and mappings under
consideration are smooth.
\par
Each function $G:A^m\to A$ can be written as
$$
     G=e_i\cdot G^i,i=0,..,n-1,
$$
where $e_i$ is a real basis of $A$. If $G$ is $A$-differentiable,
one can choose $e_i$ such that
$$
G(X^1,...,X^m)=g(x^1,...,x^m)+\sum_{|p|=1}^{n}\frac{1}{p!}\cdot
     \frac{D^pg}{dx^p}\cdot (X-x)^p,
\quad (1)
$$
where $p=(p_1,p_2,...,p_m)$, $p!=p_1!\cdot p_2!\cdot ...\cdot
p_m!,$ $|p|=p_1+p_2+...+p_m$, and
$$
(X-x)^p=(X^1-x^1)^{p_1}\cdot (X^2-x^2)^{p_2}\cdot ...\cdot
(X^m-x^m)^{p_m},
$$
$$
\frac{D^pg}{dx^p}=
\frac{\partial^{|p|}g}{(\partial x^1)^{p_1}(\partial x^2)^{p_2}...
     (\partial x^m)^{p_m}},
     X^{j}=e_i\cdot x^{j,i},i=1..n-1,
$$
$$
e_0=1\in A,(e_k)^S=0, k=1..n-1,\mbox{ S is a positive integer},
$$
i.e. $e_k$ belong to the radical $Rd(A)$ and there exists a
pseudobasis $e_l$ such that $l=1..r, r\leq n-1$,
$e_k=e_{1}^{s_{1}}\cdot..\cdot e_{r}^{s_{r}}$ (see [1],[2],[3]).
\emph{We will call this basis the standard basis.}

We will use the following proposition (see [1]).
\smallskip\par
{\bf Proposition.}[1] \emph{Suppose that  $
\dim_{ R} A = n$. Then the standard basis of $A$
satisfies: $e_0 = 1\in  A $, $ (e_1)^S = 0 $,
 $(e_k)^S = 0 $, and $ k=2..n-1 $, where $ e_k $ is a
real basis of an ideal in $A$.}
\smallskip\par
Recall also that on a manifold $M$ over $A$ carries a foliation
$\mathcal{F}$ called the \emph{canonical foliation}.

\medskip

\textbf{1. Main theorem}
\smallskip
\par
{\bf Theorem 1.} \emph{Let $M$ be a manifold over $A$,
$\dim_RA=n$, $\dim_AM=m$. Let $g$ be the real part of an
$A$-differentiable function on $M$.}

\emph{Then for any leaf $L$ of $\mathcal{F}$ there exists a leaf
$K$ of $\mathcal{F}$ such that $K\subset \overline{L}$ and
$d(g)|_K=0$.}
\smallskip
\par
\emph{Proof.} By definition, we have
$$
e_1\cdot e_k=e_i\cdot v^i_k, \mbox{ where } i,k=2..n-1,v^i_k\in
R.
$$
Evidently, $(e_1)^S=0$ implies $e_1\cdot e_1=e_i\cdot v^i$. Then
 from (1) it follows that $G(X^1,..,X^m)=g(x^1,..x^m)+e_i\cdot
G^i=g+e_i\cdot g^i$. Here $g^i$ are the coordinates of $G$ with
respect to a basis of $A$.

Let $\mathcal{A}$ be the maximal $A$-differentiable atlas. In the
general case, with respect to a chart of $\mathcal{A}$, we have
$$
g^1=\frac{\partial g}{\partial x^{j}}\cdot x^{j,1}+b(x^1,..,x^m),
\quad j=1..m.
$$

For the following reasoning it is important that
$g$ is a basic function of the canonical foliation on $M$
and $b(x^1,..,x^m)$ depends only on the
transversal coordinates $\{x^1,..,x^m\}$ of
$\mathcal{F}$ (see  [1],[2],[3],[4]).

Since $M$ is compact,  $\overline{L}$ is compact, too. Then $g^1$
has a minimum on $\overline{L}$ at a point $x \in K \subset
\overline{L}$, where $K$ is a leaf of $\mathcal{F}$. As is known,
$K \subset \overline{L}$ [4]. One can take a chart $(U, x^i)$ in
$\mathcal{A}$ containing $x$ such that  $x$ has coordinates
$\{x^1,..,x^m,0,...,0\}$. If $d(g)|_K \neq 0$, then there exits a
point $\tilde{x}$ in $U$ with coordinates
$\{x^1,..,x^m,x^{1,1},...,x^{m,1},0,...,0\}$
such that
$g^1(x)>g^1(\tilde{x})$. Since $\tilde{x} \in K$ and
$x$ is a minimum point of $g^1$,
we obtain a contradiction. $\Box$
\smallskip
\par
From Theorem 1 it follows  that  $g^1|_{\bar{L}}$ achieves maximum
and minimum on $\bar{L}/L$.
\smallskip
\par
{\bf Theorem 2.} \emph{Let $M$ be a compact
manifold over $A$, ($\dim_R A=n$, $\dim_A M=m$).
Suppose that $e_i$ is a standard basis $A$. Let G
be an $A$-differentiable function on $M$. Then}
\par
1) \emph{The real part $g$ of $G$ is constant};
\par
2) \emph{The  $e_1$-component of $G$ is constant on leaves of}
$\mathcal{F}$.
\smallskip
\par
\emph{Proof}. From Theorem 1 it follows  that all values of the
basic function $g$ are critical. Then the theorem assertion
follows from the Sard Theorem. $\Box$
\smallskip
\par
Let us denote by $\tilde{e}_k$ the element of
 the  standard basis of   $A$
such that $\tilde{e}_k \cdot e_l =0$, where $l=1..n-1$. This means
that $\tilde{e}_k \cdot Rd(A) = 0$. The other elements of the
standard basis of $A$ we will denote by $\breve{e}_j$.
\smallskip
\par
{\bf Corollary 1.} Let $M$ --- be a compact manifold over $A$, and
$G$ be an $A$-differentiable function on $M$. Then
$$
G=a+f^{k}\cdot\tilde{e}_k,
$$
where $a\in A$, and $f^k$ are basic functions.
\smallskip\par
\emph{Proof.} The corollary follows from Theorem 2 and the
properties of the standard basis $\{e_i\}$. $\Box$
\medskip
\par
\textbf{2. Dimensions of 1-form spaces on $M$.}
\smallskip
\par
Let $M$ be a manifold over $A$, $\Omega$ be the space of
differentiable forms on $M$. Recall that an $\omega \in \Omega$ is
said to be $A$-differentiable if $\omega$ and $d\omega$ are
$A$-linear. Let us denote by $\Omega_A$ the space of
$A$-differentiable forms on $M$. Note that
$$
\Omega_A\subset
\Omega_{A-lin}\subset A\otimes\Omega,
$$
where $\Omega_{A-lin}$ is the space of $A$-linear forms[2].

From this we get the following short exact sequence of complexes:
$$
0\to\Omega_A\to A\otimes\Omega\to
\frac{A\otimes\Omega}{\Omega_A}\to 0.
$$
The left complex $(\Omega_A,d)$ is called the complex of
$A$-differentiable forms. Then we obtain the cohomology exact
sequence
$$
0\to H^0_A\to A\otimes H^0\to H^0_{/A}\to
     H^1_A\to A\otimes H^1\to ... ,
$$
where the groups without lower indices are the de Rham cohomology
groups of $M$, $H^{*}_{A}$ are the cohomology groups of
$(\Omega_{A},d)$, and $H^{*}_{/A}$ are the cohomology groups of
the quotient complex $(\Omega_{/A},d)$.
\smallskip
\par
{\bf Lemma 1.} $\dots \to 0 \to H^1_A \to A\otimes H^1\to \dots$
is a  monomorphism.
\smallskip
\par
For every $j_0$, let
$$
\breve{\Omega}_{A_R}^{s,j_0} = \{\breve{\omega}^{j_0} \in \Omega^s
| \exists  \omega \in \Omega_A , \omega= \breve{e}_{j_1}\cdot
\breve{\omega}^{j_1}+ \breve{e}_{j_0}\cdot \breve{\omega}^{j_0}+
\breve{e}_{j_2}\cdot \breve{\omega}^{j_2}+ \tilde{e}_k \cdot
\tilde{\omega}^k\}.
$$
Note that $\breve{\Omega}_{A_R}^{s,j_0}\subset\Omega^s$.
\par
For each $j_0$, let
$$
Z\breve{\Omega}_{A_R}^{s,j_0} =
\{\breve{\omega}^{j_0} \in \Omega^s | \exists  \omega \in
(\Omega_A
\cap \ker\{d\}),
\omega=
\breve{e}_{j_1}\cdot \breve{\omega}^{j_1}+
\breve{e}_{j_0}\cdot \breve{\omega}^{j_0}+
\breve{e}_{j_2}\cdot \breve{\omega}^{j_2}+
\tilde{e}_k \cdot \tilde{\omega}^k\}
$$
Note that
$Z\breve{\Omega}_{A_R}^{s,j_0}\subset\breve{\Omega}_{A_R}^{s,j_0}$.
\smallskip
\par
{\bf Theorem 3.} \emph{Let $M$ be a compact manifold over $A$,
$\dim_R A=n$, and $\dim_A M=m$. Then the spaces
$Z\breve{\Omega}_{A_R}^{1,j}$
 have the finite dimension and}
$$
\dim\{Z\breve{\Omega}_{A_R}^{1,j}\}\leq
\dim\{ A\otimes H^1(M)\}.
$$
\smallskip
\par
\emph{Proof.} Let $\omega, \varsigma \in [\omega] \in  H^1_A $,
and
$\omega \in\Omega_A$,
$\omega=e_i\cdot\omega^i=
\breve{e}_j\cdot \breve{\omega}^j + \tilde{e}_k \cdot
\tilde{\omega}^k$,
$\varsigma \in\Omega_A$,
$\varsigma=e_i\cdot\varsigma^i= \breve{e}_j\cdot
\breve{\varsigma}^j + \tilde{e}_k \cdot \tilde{\varsigma}^k$.
From Corollary 1 it follows 1 that
$\breve{\Omega}_{A_R}^{0,j}
\cong R$.
Hence
$\breve{\omega}^j=\breve{\varsigma}^j$ and
$ \dim\{ H^1(M)\} <\infty $.

The second statement follows from Lemma 1.
$\Box$
\bigskip
\par
{\bf References.}
\par
   [1]  V.V. Vishnevskii, A.P. Shirokov, V.V. Shurygin.
     \emph{Spaces over algebras},
     Kazan University Press.1985.
\par
   [2] Shurygin, V. V. \emph{Jet bundles as manifolds
 over algebras}. (Russian) Translated in J. Soviet Math. 44 (1989),
  no. 2, 85--98. Itogi Nauki i Tekhniki, Problems of geometry,
   Vol. 19 (Russian), 3--22, Akad. Nauk SSSR, Vsesoyuz.
   Inst. Nauchn. i Tekhn. Inform., Moscow, 1987.
(MR0933562 (89d:58006)).
     \par [3] Shurygin,  V. V.\emph{ Manifolds over algebras and their
 application in the geometry of jet bundles.}
 (Russian) Uspekhi Mat. Nauk 48 (1993), no. 2(290),
  75--106; translation in Russian Math. Surveys 48 (1993),
  no. 2, 75--104  (MR1239861 (94i:58004)).
\par
     [4] Molino P. {\em Riemannian foliations}. -- Birkh\"auser, 1988.

\address{N.\,G.\,Chebotarev Research Institute of Mathematics and
Mechanics, Geometry Department, Universitetskaya,
17, Kazan, TATARSTAN: 420008, RUSSIAN FEDERATION}

\email{Tagir.Gaisin@ksu.ru}

\end{document}